\documentclass[11pt]{amsart}  
\usepackage{amsmath,amssymb,amsfonts,hyperref,times}

\title{Optimality and uniqueness of the (4,10,1/6) spherical code}

\author{Christine Bachoc} 
\address{C. Bachoc, Laboratoire A2X, Universit\'e Bordeaux I, 351,
cours de la Li\-b\'e\-ration, 33405 Talence, France}
\email{bachoc@math.u-bordeaux1.fr}

\author{Frank Vallentin} 
\address{F. Vallentin, Centrum voor Wiskunde en Informatica (CWI),
Kruislaan 413, 1098 SJ Amsterdam, The Netherlands}
\email{f.vallentin@cwi.nl}

\thanks{The second author was supported by the Netherlands
Organization for Scientific Research under grant NWO 639.032.203 and
by the Deutsche Forschungsgemeinschaft (DFG) under grant SCHU
1503/4-2.}


\date{May 7, 2008}

\newtheorem{defi}{Definition}[section]

\newtheorem{theorem}[defi]{Theorem}

\newcommand{\R}{{\mathbb{R}}}

\newcommand{\Sn}{S^{n-1}}

\newcommand{\trace}{\operatorname{trace}}

\begin{document}

\begin{abstract}
  Linear programming bounds provide an elegant method to prove optimality
  and uniqueness of an $(n,N,t)$ spherical code. However, this method
  does not apply to the parameters $(4,10,1/6)$. We use semidefinite
  programming bounds instead to show that the Petersen code, which
  consists of the midpoints of the edges of the regular simplex in
  dimension $4$, is the unique $(4,10,1/6)$ spherical code.
\end{abstract}

\maketitle

\section{Introduction}

Let $C$ be an $N$-element subset of the unit sphere $\Sn \subseteq
\R^n$. It is called an {\em $(n, N, t)$ spherical code} if every two
distinct points $(c,c')$ of $C$ have inner product $c \cdot c'$ at
most $t$. An $(n,N,t)$ spherical code is called {\em optimal} if there
is no $(n,N,t')$ spherical code with $t' < t$.

Only for a few parameters optimal spherical codes are known. The table
\cite[page 115]{CFG} lists all known cases in dimension $n = 3$. The
tables \cite[Table 9.1]{Lev} and \cite[Table 1]{CK1} list all
known cases in which optimality can be proven using linear programming
bounds.

One source of optimal spherical codes are iterated kissing
configurations coming from the $E_8$ root lattice in dimension $8$ and
the Leech lattice in dimension $24$ (see \cite{CS}). Starting from the
sphere packing defined by these lattices one fixes one sphere and
considers all spheres in the packing touching the fixed one. The
touching points, also called a {\em kissing configuration}, form
$(8,240,1/2)$ and respectively $(24,196560,1/2)$ spherical codes. Then
one views the kissing configuration as a packing in spherical geometry
and repeats this construction. One gets $(7,56,1/3)$ and respectively
$(23,4600,1/3)$ spherical codes.

More formally, one picks a point $x \in C$ from an $(n, M, 1/k)$
spherical code $C$ in which $x$ has $M'$ points $N_x \subseteq C$ with
inner product $1/k$. Then the points $(N_x - x/k)/\sqrt{1-1/k^2}$ form
an $(n-1, M', 1/(k+1))$ spherical code.

In this way one gets sequences of spherical codes with parameters
\[
(8,240,1/2), (7,56,1/3), (6,27,1/4), (5,16,1/5), (4,10,1/6), (3,6,1/7),
\]
and
\[
(24,196560,1/2), (23, 4600, 1/3), (22,891,1/4),
(21,336,1/5), (20,170,1/6).
\]

By using linear programming bounds Levenshtein \cite{Lev} proved that
every sharp (see Section~\ref{lp}) spherical code is optimal. Levenshtein's
theorem applies to all spherical codes above except to those with
parameters $(4,10,1/6)$, $(3,10,1/7)$, $(21,336,1/5)$,
$(20,170,1/6)$. In all optimal cases the spherical code is also
unique up to orthogonal transformations. This was proved for the cases
$(8,240,1/2)$, $(7,56,1/3)$, $(24,196560,1/2)$, $(23,4600,1/3)$ by
Bannai and Sloane \cite{BS} and for $(22,891,1/4)$ by
Cuy\-pers \cite{Cuy} and independently by Cohn and Kumar \cite{CK2}
(who also corrected a minor error in the $(23,4600,1/3)$ case). For
the cases $(6,27,1/4)$, $(5,16,1/5)$ see the discussion in \cite[Appendix
A]{CK1}. One should point out that optimality does not imply
uniqueness as one can see from the sharp $(q(q^3+1)/(q+1),
(q+1)(q^3+1), 1/q^2)$ spherical codes from \cite{CGS}. For some $q$
there are two different spherical codes with these parameters.

Based on massive computer experiments Cohn et al. \cite[Section
3.4]{BBCGKS} conjectured that the $(4,10,1/6)$ spherical code is
optimal and unique. As we explain in Section~\ref{petersen} this
spherical code is closely related to the Petersen graph and we call it
the \textit{Petersen code}. Whether the above spherical codes with
parameters $(21, 336, 1/5)$ and $(20, 170, 1/6)$ are optimal and
unique is currently unclear. At least in all these cases linear
programming bounds cannot be used to show optimality. A $(3,6,1/7)$
spherical code is not optimal because the vertices of the regular
octahedron form a $(3,6,0)$ spherical code which is a sharp spherical
code.

The main result of this paper is the following theorem which proves
the conjecture.

\begin{theorem}
\label{petersen thm}
The Petersen code is an optimal $(4,10,1/6)$ spherical code. Up to
orthogonal transformations it is the unique spherical code with these
parameters.
\end{theorem}

The proof is based on the semidefinite programming bounds for
spherical codes developed in \cite{BV} and \cite{BV2}. Currently this
is the only new case we know where the semidefinite programming bound
is tight and the linear programming bound is not. Another known case
seems to be $8$ points in $S^2$ which was solved by Sch\"utte and van
der Waerden in \cite{SW}. The linear programming bound gives $8.29$
whereas our numerical calculations suggest that the semidefinite
programming bound is tight.

We could not prove optimality of $(21,336,1/5)$ and $(20, 170, 1/6)$
spherical codes using semidefinite programming bounds.  For the first
case the linear programming bound equals $392$ whereas our numerical
calculations suggest that the semidefinite programming bound is
approximately $363$. However, we run into serious numerical problems
here and at the moment we cannot definitely rule out that the
semidefinite programming bound is sharp. For the second case the
linear programmig bound and the semidefinite programming bound
coincide: They both give $206.25$.

The structure of the paper is as follows: After giving some
constructions and properties of the Petersen code in
Section~\ref{petersen}, which also reveal the origin of its name, we
show in Section~\ref{lp} that one {\em cannot} prove
Theorem~\ref{petersen thm} using linear programming bounds. In
Section~\ref{sdp} we recall the semidefinite programming bounds and in
Section~\ref{optuniq} we present a proof of Theorem~\ref{petersen thm}
based on them.

\section{Constructions and properties of the Petersen code}
\label{petersen}

There are many possibilities to construct the Petersen code and we
already gave one. Here we give two more.

The next construction justifies the name ``Petersen code''. The
Petersen graph is a graph with $10$ vertices and $15$ edges. The
vertices are given by the $2$-element subsets of a $5$-element set and
they are adjacent whenever the corresponding $2$-element subsets have
empty intersection. Every point of the Petersen code corresponds to a
vertex of the Petersen graph and the inner product between two points
is $-2/3$ whenever the corresponding vertices are adjacent. The inner
product is $1/6$ whenever the corresponding vertices are not
adjacent. This defines a Gram matrix having rank $4$ which is unique
up to simultaneous permutation of rows and columns. The number of
ordered pairs in the Petersen code with inner product $-2/3$ is $30$
and those with inner product $1/6$ equals $60$.

In the Petersen graph every vertex has three neighbors, every pair of
adjacent vertices has no common neighbors and every pair of
nonadjacent vertices has exactly one common neighbor. So it is a
strongly regular graph with parameters $\nu = 10$, $k = 3$, $\lambda =
0$, $\mu = 1$. It is easy to see that it is uniquely defined by these
parameters. For more information about strongly regular graphs see
\cite{BI} and \cite{BCN}.

The next construction is geometric: After applying a suitable
similarity transformation the midpoints of the edges of the regular
simplex in dimension $4$ form the Petersen code. Sometimes, this
construction is also called the {\em second hypersimplex}
$\Delta(2,5)$. The second hypersimplex is the $4$-dimensional polytope
defined as the convex hull of the points $e_i + e_j$ with $1 \leq i <
j \leq 5$ where $e_i$ is the $i$-th standard unit vector in
$\R^5$. For more information about second hypersimplices see
\cite{LST}.

By \cite[Theorem 5.5]{DGS} the Petersen code forms a spherical
$2$-design: A spherical code $C \subseteq \Sn$ forms a {\em spherical
  $M$-design} if for every polynomial function $f : \R^n \to \R$ of
degree at most $M$, the average over $C$ equals the average over the
sphere $\Sn$.

\section{Linear programming bounds}
\label{lp}

Linear programming bounds provide an elegant method to prove optimality
and uniqueness of an $(n,N,t)$ spherical code. In particular a theorem
of Levenshtein \cite[Theorem 1.2]{Lev}, which covers many cases in a
unified way is based on them. Before we prove that linear programming
bounds cannot prove the optimality of the $(4,10,1/6)$ spherical code
we briefly review the underlying notions (see also e.g. \cite[Theorem
4.3]{DGS}, \cite{KL}, \cite[Chapter 9]{CS}, \cite[Theorem 2.1]{BV}).

The positivity property of the Gegenbauer polynomials $C^{n/2-1}_k$
(see \cite[Chapter 6.4]{AAR}), which are normalized by $C^{n/2-1}_k(1)
= 1$, underlies the linear programming bounds for spherical codes in
$\Sn$: For every degree $k = 0, 1, \ldots$ and every finite subset $C$
of $\Sn$ we have
\begin{equation}
\label{pos prop gegenbauer}
\sum_{(c,c') \in C^2} C^{n/2-1}_k(c \cdot c') \geq 0.
\end{equation}
One formulation of the linear programming bounds is as follows.

\begin{theorem}
\label{lpbounds}
Let $F(x)$ be a polynomial with expansion
\begin{equation}
F(x) = \sum_{k = 0}^d f_k C^{n/2-1}_k(x)
\end{equation}
in terms of Gegenbauer polynomials $C^{n/2-1}_k$. Suppose that
\begin{enumerate}
\item[(a)] all coefficients $f_k$ are nonnegative,
\item[(b)] $f_0 > 0$,
\item[(c)] $F(x) \leq 0$ for all $x \in [-1,t]$. 
\end{enumerate}
Then an $(n,N,t)$ spherical code satisfies
\begin{equation}
\label{lpboundineq}
N \leq \frac{F(1)}{f_0}.
\end{equation}
\end{theorem}

\begin{proof}
For an $(n,N,t)$ spherical code $C$ we have the inequalities
\begin{equation}
\label{lp proof}
N F(1) \geq \sum_{(c,c) \in C^2} F(1) + \sum_{\stackrel{(c,c') \in C^2}{c \neq
c'}} F(c \cdot c') = \sum_{(c,c') \in C^2} F(c \cdot c') \geq N^2
f_0.
\end{equation}
where the first inequality is due to (c) and the second due to (a) and
the positivity property \eqref{pos prop gegenbauer}. This together
with (b) implies \eqref{lpboundineq}.
\end{proof}

If there exists an $(n, N, t)$ spherical code $C$ so that $N = \lfloor
F(1)/f_0 \rfloor$ in \eqref{lpboundineq}, then, of course, $C$ is a
maximal $(n, N, t)$ spherical code, i.e.\ $N$ is the maximal number of
points one can place on the sphere $S^{n-1}$ so that distinct points
have inner product at most $t$. If furthermore $N = F(1)/f_0$, then
$C$ is an optimal $(n, N, t)$ spherical code. This can be seen as
follows. If \eqref{lpboundineq} is tight it follows from the proof
that for an $(n,N,t)$ spherical code $C$ one has $F(c \cdot c') = 0$
for distinct $c,c' \in C$. Suppose $C'$ is an $(n,N,t')$ spherical
code with $t' < t$. Then, $F(c \cdot c') = 0$ for all distinct $c,c'
\in C'$. Now we perturb $C'$ continuously to another $(n,N,t'')$
spherical code $C''$ with $t' < t'' < t$. Still we would have that $c
\cdot c'$ is a root of the polynomial $F$ for all distinct $c,c' \in
C''$ yielding a contradiction.

Levenshtein's theorem says that for every sharp spherical code there
is a polynomial satisfying the assumptions of Theorem~\ref{lpbounds}
for which \eqref{lpboundineq} is tight. A spherical code $C$ is called
{\em sharp} if it is a spherical $M$-design and the number $m$ of
different inner products between distinct points satisfies $M \geq 2m
- 1$. The Petersen code is a spherical $2$-design and there are $2$
different inner products between distinct points. Thus, Levenshtein's
theorem does not apply to it. 

Now we show that it is not possible to prove the optimality of the
Petersen code with help of Theorem~\ref{lpbounds}. Suppose that the
polynomial $F(x) = 1 + \sum_{k = 1}^d f_k C^1_k(x)$ satisfies $f_k
\geq 0$ for $k = 1, \ldots, d$ and $F(x) \leq 0$ for all $x \in
[-1,1/6]$. If $F$ would prove that the Petersen code is optimal, then
the inequalities in \eqref{lp proof} are equalities, so we would have
that
\begin{equation}
\label{first cond}
10 = F(1) = 1 + \sum_{k=1}^d f_k,
\end{equation}
and that
\begin{equation}
\label{second cond}
0 = F(-2/3) = F(1/6),
\end{equation}
and furthermore that for all $k$ with $f_k > 0$
\begin{equation}
\label{pairs eq}
0 = \sum_{(c,c') \in C^2} C^1_k(c \cdot c') = 10 +
30C^1_k(-2/3) + 60C^1_k(1/6).
\end{equation}

We shall show that \eqref{pairs eq} only holds for $k = 1$ and $k =
2$: By \cite[(6.4.11)]{AAR} we have the following expression
\begin{equation}
C^1_k(\cos \theta) = \frac{1}{k+1} \sum_{j = 0}^k \cos((k-2j)\theta).
\end{equation}
Hence,
\begin{equation}
\lim_{k \to \infty} C^1_k(-2/3) = \lim_{k \to \infty} C^1_k(1/6) = 0,
\end{equation}
so that for sufficiently large $k$, \eqref{pairs eq} cannot hold
true. Checking the remaining cases it follows that \eqref{pairs eq} is
only valid for $k = 1, 2$.  Hence, $F$ is of degree $2$, but then $F$
cannot satisfy the conditions \eqref{first cond} and \eqref{second
cond} and $F(x) \leq 0$ for $x \in [-1,1/6]$.

This argument gives rather pessimistic estimates. In fact, numerical
computations suggest that for all $d \geq 3$ the optimal polynomial is
\begin{equation}
F(x) = 1 + \frac{2270}{680} x + \frac{2775}{680} \big(\frac{4}{3}x^2 -
\frac{1}{3}\big) + \frac{1500}{680} \big(2x^3 - x\big),
\end{equation}
and so the best upper bound one can probably prove using
Theorem~\ref{lpbounds} is $10.625$.  We checked this for all $d
\leq 40$ by computer.

\section{Semidefinite programming bounds}
\label{sdp}

As we have seen above the positivity property of the polynomials
$C^{n/2-1}_k$ plays a crucial role for the linear programming
bounds. For the semidefinite programming bounds this is replaced by
the positivity property of the matrices $S^n_k$. From \cite{BV} we
recall the matrices $S^n_k$ and their positivity property. First we
define the entry $(i,j)$ with $i, j \geq 0$ of the (infinite) matrix
$Y^n_k$ containing polynomials in $x,y,z$ by
\begin{equation}
\label{ynk}
\begin{split}
\big(Y^n_k\big)_{i,j}(x,y,z) & = x^i y^j \;\cdot\\ &((1-x^2)(1-y^2))^{k/2}
C^{n/2-3/2}_k\left(\frac{z-xy}{\sqrt{(1-x^2)(1-y^2)}}\right),
\end{split}
\end{equation}
and then we get $S^n_k$ by symmetrization:
\begin{equation}
S^n_k = \frac{1}{6}\sum_{\sigma} \sigma Y^n_k,
\end{equation}
where $\sigma$ runs through all permutations of the variables $x,y,z$ which acts on the matrix coefficients in the obvious way.
The matrices $S^n_k$ satisfy the positivity property:
\begin{equation}
\text{for all finite $C \subseteq \Sn$,} \sum_{(c,c',c'') \in C^3}
S^n_k(c \cdot c', c\cdot c'', c' \cdot c'') \succeq 0,
\end{equation}
where ``$\succeq 0$'' stands for ``is positive semidefinite'' where we
mean that every finite minor is positive semidefinite. Note that the
difference between \eqref{ynk} and the original \cite[(11)]{BV} is due
to a change of basis which does not affect the positivity property.

The interval $[-1,t]$ of the linear programming bounds is supplemented
by the domain
\begin{equation}
D = \{(x,y,z) : -1 \leq x, y, z \leq t, 1 + 2xyz - x^2 - y^2 -
z^2 \geq 0\}.
\end{equation}

We need some more notation. The space of (finite) symmetric matrices
is a Euclidean space with inner product $\langle F, G \rangle =
\trace(FG)$. The cone of positive semidefinite matrices is self dual,
i.e.\ one has $\langle F, G\rangle \geq 0$ for all positive
semidefinite $G$ if and only if $F$ is positive semidefinite. If $F$
is a symmetric matrix with $m$ rows and $m$ columns, then we interpret
$\langle F, S^n_k\rangle$ as the inner product of $F$ with the
principal minor of $S^n_k$ of appropriate size.

Now we can state the semidefinite programming bounds. The following
polynomial formulation can be deduced from \cite[Theorem 4.2]{BV}. We
provide an independent proof which has the additional feature that it
gives information in the case when the theorem provides tight results.

\begin{theorem}
\label{sdpbounds}
Let $F(x,y,z)$ be a symmetric polynomial with expansion
\begin{equation}
\label{snk expansion}
F(x,y,z) = \sum_{k = 0}^d \langle F_k, S^n_k \rangle,
\end{equation}
in terms of the matrices $S^n_k$. Suppose that 
\begin{enumerate}
\item[(a)] all $F_k$ are positive semidefinite
\item[(b)] $F_0 - f_0 E_0 \succeq 0$ for some $f_0 > 0$ ($E_0$ is the
matrix whose only nonzero entry is the top left corner which contains
$1$),
\item[(c)] $F(x,y,z) \leq 0$ for all $(x,y,z) \in D$,
\item[(d)] $F(x,x,1) \leq B$ for all $x \in [-1,t]$.
\end{enumerate}
Then an $(n,N,t)$ spherical code satisfies
\begin{equation}
\label{bound}
N \leq \frac{3B + \sqrt{9B^2 + 4f_0(F(1,1,1) - 3B)}}{2f_0}.
\end{equation}
\end{theorem}

\begin{proof}
Let $C$ be an $(n,N,t)$ spherical code. Define
\begin{equation}
S = \sum_{(c,c',c'') \in C^3} F(c \cdot c', c \cdot c'', c' \cdot c'').
\end{equation}

Split this sum into three parts according to the indices $C_1, C_2,
C_3 \subseteq C^3$ where $C_i$ contains all triples with $i$ pairwise
different elements. The contribution of $C_1$ to $S$ is $NF(1,1,1)$,
the one of $C_2$ at most $3N(N-1)B$ and the one of $C_3$ is at most
zero. Together,
\begin{equation}
S \leq N F(1,1,1) + 3N(N-1)B.
\end{equation}

On the other hand,
\begin{eqnarray}
S & = & \sum_{k = 0}^d \langle F_k, \sum_{(c,c',c'') \in C^3} S^n_k(c \cdot c', c \cdot c'', c' \cdot c'')\\
& \geq & \langle f_0 E_0, \sum_{(c,c',c'') \in C^3} S^n_k(c \cdot c', c \cdot c'', c' \cdot c'') \rangle\\
& = & N^3 f_0,
\end{eqnarray}
yielding the statement of the theorem.
\end{proof}

A few remarks about the theorem and its proof are in order.

\medskip

If the bound \eqref{bound} is tight, then all inequalities in the
proof must be equalities. In particular, the univariate polynomial
$F(x,x,1) - B$ has roots at the inner products $c \cdot c'$ for
distinct $c, c' \in C$. So we can argue in the same way as in the case
of the linear programming bounds that tightness implies
optimality.

\medskip

If the bound \eqref{bound} is tight, we have the following
identities: Let $C$ be an $(n,N,t)$ spherical code with
\begin{equation}
\begin{split}
D(C) & = \{(c \cdot c', c \cdot c'', c' \cdot c'') : (c,c', c'') \in C^3\},\\
I(C) & = \{c \cdot c' : (c,c') \in C^2, c \neq c'\}.
\end{split}
\end{equation}
Let $F$ be a polynomial satisfying the hypothesis of
Theorem~\ref{sdpbounds} with constants $B$ and $f_0$ and proving the 
tight bound $(3B + \sqrt{9B^2 + 4f_0(F(1,1,1) - 3B)})/2f_0$. Then
\begin{enumerate}
\item[(i)] $N^2f_0 - F(1,1,1) - 3(N-1)B= 0$,
\item[(ii)] $F(x,y,z) = 0$ for all $(x,y,z) \in D(C)$,
\item[(iii)] $F(x,x,1) = B$ for all $x \in I(C)$,
\item[(iv)] $\langle F_k, \sum_{(c,c',c'') \in C^3} S^n_k(c \cdot c', c \cdot c'', c' \cdot c'')\rangle = 0$ for all $k = 1, \ldots, d$,
\item[(v)] $\langle F_0, \sum_{(c,c',c'') \in C^3} S^n_0(c \cdot c', c \cdot c'', c' \cdot c'')\rangle = N^3 f_0$.
\end{enumerate}

\medskip

Semidefinite programming bounds are at least as strong as linear
programming bounds: If $G = \sum_{k = 0}^d g_k C^{n/2-1}_k(x)$ is a
polynomial which satisfies the hypothesis of Theorem~\ref{lpbounds},
then the polynomial $F(x,y,z) = (G(x) + G(y) + G(z))/3$ satisfies the
hypothesis of Theorem~\ref{sdpbounds} with $B = G(1)/3$ and $f_0 =
g_0$. This is because one sets $F_0 = g_0 E_0$ and from
\cite[Proposition 3.5]{BV} it follows that one can express $G$ with
semidefinite matrix coefficients.

\medskip

From \cite[Lemma 4.1]{BV2} it follows that one can express every symmetric
polynomial in the form \eqref{snk expansion}. However, this expansion is not
unique, e.g.\ 
\begin{equation}
\begin{split}
x + y + z & = 
\langle
\left(
\begin{smallmatrix}
0 & 3/2\\
3/2 & 0
\end{smallmatrix}
\right),
S^n_0
\rangle +
\langle 
\left(
\begin{smallmatrix}
0
\end{smallmatrix}
\right),
S^n_1
\rangle\\
&  =
\langle
\left(
\begin{smallmatrix}
0 & 0\\
0 & 3\\
\end{smallmatrix}
\right),
S^n_0
\rangle +
\langle 
\left(
\begin{smallmatrix}
3
\end{smallmatrix}
\right),
S^n_1
\rangle,
\end{split}
\end{equation}
where only the second expansion involves semidefinite matrices and
where
\begin{equation}
S^n_0 =
\left( 
\begin{smallmatrix}
1 & (x + y + z)/3\\
(x + y + z)/3 & (xy + xz + yz)/3
\end{smallmatrix}
\right),\;
S^n_1 = 
\left(
\begin{smallmatrix}
(x + y + z)/3 - (xy + xz + yz)/3
\end{smallmatrix}
\right).
\end{equation}

\section{Proof of optimality and uniqueness}
\label{optuniq}

In this section we prove Theorem~\ref{petersen thm} with the help of
Theorem~\ref{sdpbounds}. Although we can present a proof which one can
verify essentially without using computer we relied heavily on
computer assistance to find it.

To show that the Petersen code is the unique $(4,10,1/6)$ spherical
code we use the matrices $F_0 \in \R^{4 \times 4}$, $F_1 \in \R^{3
\times 3}$, $F_2 \in \R^{1 \times 1}$ given by
\begin{equation}
\begin{split}
F_0 & =
\begin{pmatrix}
2882/3 & 114 & -2500 & 0\\
114 & 324 & 216 & 0\\
 -2500 & 216 & 8716 & 1296\\
 0 & 0 & 1296 & 11664
\end{pmatrix},\\
F_1 & = 
\begin{pmatrix}
0 & 0 & 0\\
0 & 3588 & -4536\\
0 & -4536 & 11664
\end{pmatrix},
F_2 =
\begin{pmatrix}
2000
\end{pmatrix}.
\end{split}
\end{equation}
Let
\begin{equation}
m_{ijk} = \frac{1}{6}\sum_{\sigma} \sigma(x^iy^jz^k),\quad 0 \leq i \leq j \leq k,
\end{equation}
where $\sigma$ runs through all permutations of the variables $x,y,z$,
be the polynomial which one gets by symmetrizing $x^iy^jz^k$. Then,
\begin{equation}
\begin{split}
F(x,y,z) = &
11664m_{320} + 11664m_{221} + 
7128m_{220} - 9072m_{211}\\
& +432m_{210} - 2412m_{111} +
324m_{110} + 228m_{100} - 118/3,
\end{split}
\end{equation}
and
\begin{equation}
F(x,x,1) - B = \frac{1}{3888}\Big(x+\frac{2}{3}\Big)^2\Big(x-\frac{1}{6}\Big)\Big(x^2+\frac{4}{9}x+\frac{20}{27}\Big).
\end{equation}
It is a straight forward computation that $F$ satisfies the condition
of Theorem~\ref{sdpbounds} with $F(1,1,1) = 59750/3$, $B = 250$, $f_0
= 800/3$ so that it shows $N \leq 10$ for a $(4,N,1/6)$ spherical
code. This finishes the proof of the optimality.

Before showing uniqueness, let us describe how we derived $F_0$,
$F_1$, $F_2$.  We have
\begin{equation}
\begin{split}
S^4_0(x,y,z) & = \begin{pmatrix}
1 & m_{100} & m_{200} & m_{300}\\
m_{100} & m_{110} & m_{210} & m_{310}\\
m_{200} & m_{210} & m_{220} & m_{320}\\
m_{300} & m_{310} & m_{320} & m_{330}
\end{pmatrix},\\
S^4_1(x,y,z) & = \begin{pmatrix}
m_{100}-m_{110} & m_{110}-m_{210} & m_{210}-m_{310}\\
m_{110}-m_{210} & m_{111}-m_{220} & m_{211}-m_{320}\\
m_{210}-m_{310} & m_{211}-m_{320} & m_{221}-m_{330}
\end{pmatrix},\\
S^4_2(x,y,z) & = \begin{pmatrix}
-\frac{1}{2} + \frac{5}{2}m_{200} - 3m_{111} + m_{220}.
\end{pmatrix}
\end{split}
\end{equation}
Then $0 = \sum_{k = 0}^2 \langle K_{i,k}, S^4_k \rangle$ for
\begin{equation}
\begin{split}
K_{1,0} & = \begin{pmatrix}
0 & -\frac{1}{2} & 0 & 0\\
-\frac{1}{2} & 1 & 0 & 0\\
0 & 0 & 0 & 0\\
0 & 0 & 0 & 0
\end{pmatrix},
K_{1,1} = \begin{pmatrix}
1 & 0 & 0\\
0 & 0 & 0\\
0 & 0 & 0
\end{pmatrix},
K_{1,2} = \begin{pmatrix}
0
\end{pmatrix},
\end{split}
\end{equation}

\begin{equation*}
\begin{split}
K_{2,0} & = \begin{pmatrix}
\frac{1}{2} & 0 & -\frac{5}{4} & 0\\
0 & 0 & 0 & 0\\
-\frac{5}{4} & 0 & 2 & 0\\
0 & 0 & 0 & 0
\end{pmatrix},
K_{2,1} = \begin{pmatrix}
0 & 0 & 0\\
0 & 3 & 0\\
0 & 0 & 0
\end{pmatrix},
K_{2,2} = \begin{pmatrix}
1
\end{pmatrix}\\
K_{3,0} & = \begin{pmatrix}
0 & 0 & 0 & 0\\
0 & -1 & \frac{1}{2} & 0\\
0 & \frac{1}{2} & 0 & 0\\
0 & 0 & 0 & 0
\end{pmatrix},
K_{3,1} = \begin{pmatrix}
0 & \frac{1}{2} & 0\\
\frac{1}{2} & 0 & 0\\
0 & 0 & 0
\end{pmatrix},
K_{3,2} = \begin{pmatrix}
0
\end{pmatrix},\\
K_{4,0} & = \begin{pmatrix}
0 & 0 & 0 & 0\\
0 & 0 & -\frac{1}{2} & \frac{1}{2}\\
0 & -\frac{1}{2} & 0 & 0\\
0 & \frac{1}{2} & 0 & 0
\end{pmatrix},
K_{4,1} = \begin{pmatrix}
1 & 0 & \frac{1}{2}\\
0 & 0 & 0\\
\frac{1}{2} & 0 & 0
\end{pmatrix},
K_{4,2} = \begin{pmatrix}
0
\end{pmatrix},
\end{split}
\end{equation*}
i.e.\ the matrices $K_{i,k}$ form a basis of the kernel of the linear
map which assigns symmetric polynomials to the matrix
coefficients. From the discussion following the proof of
Theorem~\ref{sdpbounds} we know that the matrix entries have to
satisfy the equalities (i)--(v) where
\begin{equation}
\label{snks}
\begin{split}
\sum_{(c,c',c'') \in C^3} S^4_0(c \cdot c', c \cdot c'', c' \cdot c'') & = 
\begin{pmatrix}
1000 & 0 & 250 & \frac{125}{9}\\
0 & 0 & 0 & 0\\
250 & 0 & \frac{125}{2} & \frac{125}{36}\\
\frac{125}{9} & 0 & \frac{125}{36} & \frac{125}{648}
\end{pmatrix},\\
\sum_{(c,c',c'') \in C^3} S^4_1(c \cdot c', c \cdot c'', c' \cdot c'') & = 
\begin{pmatrix}
0 & 0 & 0\\
0 & 0 & 0\\
0 & 0 & 0
\end{pmatrix},\\
\\
\sum_{(c,c',c'') \in C^3} S^4_2(c \cdot c', c \cdot c'', c' \cdot c'') & = \begin{pmatrix} 0 \end{pmatrix}.
\end{split}
\end{equation}
We restrict our search to polynomials $F$ satisfying
\begin{equation}
\frac{\partial F}{\partial x}\Big(-\frac{2}{3}, -\frac{2}{3}, \frac{1}{6}\Big) = 0,
\frac{\partial F}{\partial x}\Big(-\frac{2}{3}, \frac{1}{6}, \frac{1}{6}\Big) = 0,
\frac{\partial F}{\partial x}\Big(-\frac{2}{3}, -\frac{2}{3}, 1\Big) = 0.
\end{equation}
Furthermore, we restrict our search to those polynomials lying in the
subspace of dimension $9$ spanned by
\begin{equation}
m_{320}, m_{221}, m_{220}, m_{211}, m_{210}, m_{111}, m_{110}, m_{100}, 1.
\end{equation}
The one dimensional affine subspace 
\begin{equation}
\begin{split}
F_{\gamma}(x,y,z) & = 
\big(11664m_{320}+9720m_{220}-1296m_{210}-6480m_{111}\\
& +2268m_{110} -108m_{100}-18\big) + \gamma\big(34992m_{221}-7776m_{220}\\
& -27216m_{211}+ 5184m_{210}+12204m_{111} - 5832m_{110}\\
& + 1008m_{100} - 64\big),\quad \gamma \in \R,
\end{split}
\end{equation}
satisfies all these linear equalities. We have
\begin{equation}
F_{\gamma}(x,y,z) = \sum_{k=0}^2 \langle A_k, S^4_k \rangle + \gamma
\langle B_k, S^4_k \rangle
\end{equation}
with
\begin{equation}
\begin{split}
A_0 = \begin{pmatrix}
-18 & -54 & 0 & 0\\
-54 & 2268 & -648 & 0\\
0 & -648 & 3240 & 5832\\
0 & 0 & 5832 & 0
\end{pmatrix},\\
A_1 = \begin{pmatrix}
0 & 0 & 0\\
0 & -6480 & 0\\
0 & 0 & 0
\end{pmatrix},
A_2 = \begin{pmatrix}
0
\end{pmatrix},
\end{split}
\end{equation}
and
\begin{equation}
\begin{split}
B_0 = \begin{pmatrix}
-64 & 504 & 0 & 0\\
 504 & -5832 & 2592 & 0\\
0 & 2592 & 4428 & -13608\\
 0 & 0 & -13608 & 34992
\end{pmatrix},\\
B_1 = \begin{pmatrix}
0 & 0 & 0\\
0 & 12204 & -13608\\
0 & -13608 & 34992 
\end{pmatrix},
B_2 = \begin{pmatrix}
0
\end{pmatrix}.
\end{split}
\end{equation}

In this affine subspace we want to find a polynomial which satisfies
the inequalities (c) and (d) from Theorem~\ref{sdpbounds} and which at
the same time has a representation of the form \eqref{snk expansion}
with positive semidefinite matrices $F_k$. Hence, we are left with the
problem of finding a matrix in the intersection of an affine subspace
with the cone of positive semidefinite matrices which is a basic task
in semidefinite programming. Since this problem is not known to be in
$\mathrm{NP}$ --- in fact it is the major open problem in the theory
of semidefinite programming --- it is not a priori clear that a
solution exists which one can nicely describe.

We solved these two semidefinite programming problems separately and
we used the numerical software \texttt{csdp} \cite{B} for this task:
If $0.28 \lessapprox \gamma \lessapprox 0.68$, then $F_\gamma$ satisfies (c). If
$0.18 \lessapprox \gamma \lessapprox 0.38$, then $F_\gamma$ has a representation of
the form \eqref{snk expansion} with positive semidefinite matrices. We
make the Ansatz $\gamma = \frac{1}{3}$ and try to find a nice
representation. For this we solve the semidefinite feasibility problem
\begin{equation}
A_k + \frac{1}{3} B_k + \beta_1 K_{1,k} + \beta_2 K_{2,k} + \beta_3
K_{3,k} + \beta_4 K_{4,k} \succeq 0, \quad k = 0, 1, 2,
\end{equation}
which luckily happens to have the solution $\beta_1 = \beta_3 =
\beta_4 = 0$ and $\beta_2 = 2000$.

To show uniqueness we first derive the three points distance
distribution $\alpha$ of a $(4,10,1/6)$ spherical code $C$ which is
defined by
\begin{equation}
\alpha(x,y,z) = \frac{1}{|C|}|\{(c,c',c'') \in C^3 : c \cdot c' = x, c \cdot c'' = y, c' \cdot c'' = z\}|.
\end{equation}
Since $-2/3$ and $1/6$ are the only roots of the polynomial $F(x,x,1)
- B$, these are the only inner products which can occur among distinct
points in $C$. This enables us to use (iv) and (v) together with the
relations
\begin{equation}
\begin{array}{l}
\alpha(x,y,z) = \alpha(\sigma(x,y,z)),\;\text{for all permutations $\sigma$ of $x,y,z$},\\
 \alpha(1,1,1) = 1,\\
 \sum_{(x,y,z) \in D} \alpha(x,y,z) = 100,\\
\sum_{x \in [-1,1]} \alpha(x,x,1) = 10,
\end{array}
\end{equation}
to determine $\alpha$ by solving a system of linear equations: It is
\begin{equation}
\begin{array}{ll}
\alpha(-2/3, -2/3, 1/6) = 6,&
\alpha(-2/3, -2/3, 1) = 3,\\
\alpha(-2/3, 1/6, 1/6) = 12,&
\alpha(1/6, 1/6, 1/6) = 18,\\
\alpha(1/6, 1/6, 1) = 6,&
\alpha(1, 1, 1) = 1.
\end{array}
\end{equation}
Now by \cite[Theorem 5.5]{DGS} $C$ is a spherical $2$-design. By
\cite[Theorem 7.4]{DGS} it carries a 2-class association scheme whose
valencies and intersection numbers are uniquely determined. In fact it
is a strongly regular graph with parameters $\nu = 10$, $k = 3$,
$\lambda = 0$, $\mu = 1$. This uniquely defines the Petersen graph
which finishes the proof of the uniqueness.

\section{Acknowledgments}

We thank Henry Cohn for communicating this problem at the Oberwolfach
seminar ``Sphere Packings: Exceptional Geometric Structures and
Connections to other Fields'' in November 2005 and for further helpful
discussions. We thank the two anonymous referees for useful
suggestions.

\end{document}